\documentclass[11pt]{amsart}  
\usepackage{vmargin, amssymb, pb-diagram}

\swapnumbers 
\newtheorem{theorem}{Theorem}[section]

\newtheorem{proposition}[theorem]{Proposition}

\newtheorem{corollary}[theorem]{Corollary}

\theoremstyle{definition}
\newtheorem{definition}[theorem]{Definition}

\theoremstyle{remark}

\newcommand{\noproof}{\hfill\qedsymbol}

\numberwithin{equation}{section}

\makeatletter
\def\romenumi{%
  \def\theenumi{\roman{enumi}}%
  \def\p@enumi{\theenumi}%
  \def\labelenumi{(\@roman\c@enumi)}}
\makeatother

\newif\ifShowLabels
\ShowLabelsfalse
\newdimen\mgheight
\def\marginnotes#1{%
                 \leavevmode\vadjust{\setbox0=\hbox{{\tt
                               \quad\quad {\small\textup #1}}}%
                 \mgheight=\ht0
                 \advance\mgheight by \dp0
                 \advance\mgheight by \lineskip
                 \kern -\mgheight \vbox to
                 \mgheight{\rightline{\rlap{\box0}} \vss}}}
\def\tlabel#1{\label{#1}
                \ifShowLabels \marginnotes{#1} \fi}

\setpapersize{USletter}
\setmarginsrb{25mm}{25mm}{25mm}{25mm}{12pt}{15mm}{12pt}{20mm} 

\DeclareMathOperator{\Spec}{Spec}
\DeclareMathOperator{\Prim}{Prim}
\DeclareMathOperator{\Max}{Max}
\DeclareMathOperator{\co}{co}
\DeclareMathOperator{\ann}{ann}

\renewcommand{\d}{\downarrow}
\renewcommand{\u}{\uparrow}

\begin{document}

\title{Primitive ideals in Hopf algebra extensions}

\author{Mark~C.~Wilson}  
\address{Department of Mathematics, University of Auckland, Private Bag 92019
 Auckland, New Zealand} 
\email{wilson@math.auckland.ac.nz}


\date{\today}

\begin{abstract} Let $H$ be a finite-dimensional Hopf algebra. We study  the behaviour of
primitive and maximal ideals in certain types of ring extensions determined by $H$. The main focus
is on the class of faithfully flat Galois extensions, which includes  includes
smash and crossed products. It is shown how analogous results can be obtained for the larger class
of extensions  possessing a total integral, which includes extensions $A^H\subseteq A$ when $H$ is
semisimple.  We use Passman's ``primitivity machine" to reduce the whole theory of Krull
relations for prime ideals to the case of  primitive ideals. The concept of strongly
semiprimitive Hopf  algebra is introduced and investigated. Several examples and open problems
are discussed.

\end{abstract}  
\maketitle 
\section{Introduction}
\tlabel{sec:intro}

A series of papers by Schneider
\cite{schneider;principal-homogeneous, schneider;representation-hopf, schneider;normal-basis} has
made clear the importance of the class
 of faithfully flat Hopf Galois extensions. As well as including smash and
 crossed products, such extensions arise in the study of algebraic groups 
\cite{schneider;principal-homogeneous}. A large class of interesting examples is provided by the
following  situations. Let $A$ be a Hopf algebra with a  normal Hopf subalgebra $R$ such that  
either (i) $A$ is pointed or the coradical of $A$ is cocommutative or  (ii) $R$ is central and
of finite index, and $R$ is a noetherian ring.
 Then $A$ is a faithfully  flat $H$-Galois extension of $R$, where $H=A/AR^+$ and 
$R^+$ denotes the augmentation ideal of $R$. For more examples, see 
\cite{montgomery-schneider;prime} and \cite{schneider;representation-hopf};  the former contains
an explicit quantum group example.

The abstract viewpoint of faithfully flat Hopf Galois extensions 
 makes it possible to give unified proofs of seemingly  disparate results. 
There are several technical advantages to working within the larger class of
 extensions, even if one is only interested in crossed  products.
 One of these is transitivity, as  defined in
Section~\ref{sec:basics}. Another is stability  under various natural operations.
 Thus even for connected Hopf algebras, for which every faithfully flat  Hopf
 Galois extension is in fact a crossed product \cite{bell;comodule}, there is 
value in this approach. 

 In \cite{montgomery-schneider;prime}, Montgomery and Schneider made a thorough
 examination of the  behaviour of prime ideals in faithfully flat Hopf Galois
 extensions, with particular
 attention to the Krull relations. The present paper contains an analogous development for
primitive ideals,  and discusses the relation between the two theories. In addition we rework
the presentation of some fundamental topics in 
\cite{montgomery-schneider;prime}, and give a more leisurely  treatment of others, in order to
allow future work to proceed more  smoothly. We hope to stimulate further work on this subject,
and so  several conjectures, most of them not new, are explicitly stated.

The outline of the paper is as follows. In Section~\ref{sec:basics} we recall the main
definitions and general results from  previous papers. Much of this material is contained in
\cite{montgomery-schneider;prime}. The chief novelty lies in our presentation, including the
systematic use of commutative diagrams to make explicit various reductions which simplify much
of the later work. Results and notation from this section are used heavily throughout the paper.

Section~\ref{sec:lying} concerns Krull relations in
$H$-extensions. We give a more leisurely treatment of some topics covered in
\cite{montgomery-schneider;prime}. Furthermore we introduce what we call internal Krull
relations, which seem to be technically useful. 

If all the Krull relations hold, then many properties of maximal ideals of $R$ yield
corresponding results for $A$, and vice versa. In view of the  difficulty of verifying all Krull
relations, however, it may be useful to  consider maximal ideals separately.
Section~\ref{sec:maximal} deals in some detail with basic questions about maximal and
$H$-maximal ideals.

Section~\ref{sec:modules} deals with modules and annihilators, in preparation for a discussion of
primitive ideals. We make explicit certain module-theoretic properties (the finite and semisimple
induction and restriction properties) which have proved crucial in group and enveloping algebra
arguments. We treat primitive and $H$-primitive ideals in Section~\ref{sec:primitive}. The aim
(largely fulfilled) is to prove analogues for $\Prim$ of all results about $\Spec$ from 
\cite{montgomery-schneider;prime}. In addition, we use Passman's ``primitivity machine" to
relate  results on the two spectra, and show thereby that the theory may be founded on $\Prim$.

A brief Section~\ref{sec:totint} discusses the relaxation of our standing hypotheses to the
larger class of  $H$-extensions with a total integral. This is  equivalent to studying
$H$-module algebras with surjective trace map, a condition which is always satisfied when $H$ is
semisimple. Because this procedure has been done many times before, and is extensively discussed
in
\cite{montgomery-schneider;prime}, we omit many details. Section~\ref{sec:radicals} discusses
various obviously defined $H$-radicals and their relation to other properties. In particular we
focus on  a property of semisimple $H$, strong semiprimitivity, which relates to the module
properties mentioned above.

In Section~\ref{sec:examples},  we discuss several classes of Hopf algebras, and summarize which
properties discussed in this paper they are known to satisfy. One new result here, answering a
question from \cite{montgomery-schneider;prime}, is that every pointed 
$H$ satisfies the Krull relation $t$-LO. 
 We do not pursue specific examples of Hopf algebras or extensions in detail, as this paper's
goal is foundational. However, there should be several quantum group applications of the content
of this paper. We conclude with some open problems in  Section~\ref{sec:conc}.

\subsection*{Notation and conventions} Let $F$ be a field and let $H$ be a Hopf algebra over
$F$. We shall be concerned with $H$-extensions, that is, extensions $R\subseteq A$ of
$F$-algebras such that $A$ is an $H$-comodule algebra via $\rho:A\to A\otimes H$ and $R=A^{\co
H}$ is its algebra of coinvariants. If 
$H$ is finite-dimensional, then this is equivalent to $A$ being an $H^*$-module and 
$R$  the algebra of invariants for $H^*$. The notation $\mathcal{E}=(R,A)$ will be used for
$H$-extensions.

A general reference for all unproved assertions about Hopf algebras which cannot  be found in
\cite{montgomery-schneider;prime} is the book \cite{montgomery;hopf-action-rings}.
 
\section{Faithfully flat Galois extensions}
\tlabel{sec:basics}

In this section we consider the fundamentals of faithfully flat 
Galois extensions. The results in this section allow for great 
simplification of proofs later on, and will be used extensively.

The $H$-extension $\mathcal{E}=(R,A)$ is said to be {\em Galois} if the map 
$A\otimes_R A\to A\otimes H$ given by $x\otimes y\mapsto x\rho(y)$ is bijective. If in addition
$A$ is faithfully flat as a left $R$-module, we say that 
$\mathcal{E}$ is {\em faithfully flat $H$-Galois}. If $H$ is  finite-dimensional and we do not
wish to specify $H$, we shall say that 
$\mathcal{E}$ is a faithfully flat finite Galois extension.

To say that $H$ has a given property (normally associated with extensions) will mean that all
faithfully flat $H$-Galois extensions  have the property.

Suppose now that $(R,A)$ is faithfully flat $H$-Galois. A fundamental fact
\cite{schneider;principal-homogeneous} is  that there is an equivalence between the category of
left $R$-modules and the category of left $A$-modules which are also right $H$-comodules, given
by $V\mapsto AV$ and $W\mapsto W^{\co H}$.  We shall use a consequence of  this heavily in the
special case of ideals.

An ideal $I$ of $R$ is said to be {\em
$H$-stable} if $IA=AI$.  If the antipode of $H$ is bijective, and the extension is a crossed
product, then the  usual definition of
$H$-stable implies the above property. For general $H$ this is no longer the case, however.

The set of all $H$-stable ideals of $R$ forms a sublattice
$\mathcal{I}_H(R)$ of the lattice  $\mathcal{I}(R)$ of all ideals of
$R$.

\begin{definition} For each ideal $I$ of $R$,  the {\em $H$-core} $(I:H)$ is the largest
$H$-stable ideal of $R$ contained in $I$.
\end{definition}

The $H$-core is well-defined since the sum of $H$-stable ideals is
$H$-stable. Taking the $H$-core gives a map $\mathcal{I}(R)\to
\mathcal{I}_H(R)$ which is a lattice epimorphism. The set-theoretic kernel of this map is a
partition of $\mathcal{I}(R)$ whose associated equivalence relation we denote by
$\sim_H$. Explicitly, $Q\sim_H Q'$ if and only if $(Q:H)=(Q':H)$. Note that each equivalence
class has a unique $H$-stable member, the common
$H$-core of all elements in the class.  

The category equivalence above, when restricted to ideals, leads to a fundamental correspondence
between $H$-stable ideals of $R$ and ideals of $A$ which are $H$-subcomodules. When $H$ is
finite-dimensional, this latter set coincides with the the set of $H^*$-stable ideals of
$A$. In any case, we shall denote it by $\mathcal{I}_{H^*}(A)$. The correspondence is given by
expansion ($I\mapsto IA$) and contraction ($J\mapsto J\cap R$). 
The following result from \cite{montgomery-schneider;prime} gives the 
precise information we require.

\begin{proposition}
\tlabel{fundcorr} Let $(R,A)$ be a faithfully flat $H$-Galois extension. Then  expansion and
contraction yield maps $\Phi: \mathcal{I}_H(R)\to
\mathcal{I}_{H^*}(A)$ and  $\Psi: \mathcal{I}_{H^*}(A)\to
\mathcal{I}_H(R)$ which are mutually inverse lattice isomorphisms.
\noproof
\end{proposition}

\subsection*{Stability properties} We now discuss the stability of faithfully flat Galois
extensions under  certain operations. In the following, a property of ring extensions is 
identified with the class of all extensions possessing the property.
\begin{definition}
 Let
$\mathcal{P}$ be a property  of ring extensions.  Say that
$\mathcal{P}$ is {\em transitive} if whenever  $R\subset B\subset A$ and both $(R,B)$ and
$(B,A)$ have the property then  $(R,A)$ has the property. If $\mathcal{E}$ is an extension of
$F$-algebras, and $E$ is a subfield or extension field of $F$,  then  $\mathcal{P}$ {\em holds
over $E$} if and only if the extension formed from $\mathcal{E}$  by (respectively) restriction
or extension of scalars to $E$ satisfies 
$\mathcal{P}$. If $\mathcal{P}$ holds over all subfields and extension  fields, then we say that
$\mathcal{P}$ is {\em field-independent}.
\end{definition}

In the next result, if whenever all factors of a subnormal series for $H$  possess a given
property $\mathcal{P}$, then $H$ has $\mathcal{P}$, we say  that $\mathcal{P}$ {\em ascends via
subnormal series}.

\begin{proposition}
\tlabel{subnormal} Let $(R,A)$ be a faithfully flat $H$-Galois extension, where $H$ is 
finite-dimensional. If $K$ is a normal Hopf subalgebra of $H$ with quotient $\overline{H}$, then 
there is an intermediate subalgebra $B$ such that $(R,B)$ is faithfully flat
 $K$-Galois and $(B,A)$ is faithfully flat $\overline{H}$-Galois. Consequently, transitive
properties ascend via subnormal series.
\noproof 
\end{proposition}

\begin{proposition} 
\tlabel{stability} The property of being faithfully flat finite Galois  holds over extension
fields and subfields, and respects quotients by  stable ideals.  More precisely, let $(R,A)$ be
a faithfully flat $H$-Galois extension,  let $I$ be an $H$-stable ideal of $R$, let 
$E_1$ be a subfield of $F$ and let $E_2$ be an
 extension field of $F$. Then
\begin{enumerate}
\romenumi
\item $(R', A')$ is faithfully flat $H'$-Galois, where $'$ denotes extension of the ground field
to $E_2$
\item $(R,A)$ is faithfully flat $H$-Galois when 
$R,A,H$ are considered as $E_1$-algebras
\item $(R/I, A/AI)$ is  faithfully flat $H$-Galois, where $R/I$ is identified with its image in 
$A/AI$ under the canonical map $A\to A/AI$.
\end{enumerate}
\noproof
\end{proposition}

\subsection*{Duality}

If $H$ is finite-dimensional, we shall let $\mathcal{E}^*$ denote the {\em dual extension} 
$(A,A\# H^*)$. Note that $H$ acts on $A\#H^*$ by  acting trivially on $A$ and with the usual
action on $H^*$. Thus we can  form $A\#H^*\#H$ which is Morita-equivalent to $A$.

\begin{definition}Let $\mathcal{P}$ be a property of faithfully flat finite Hopf Galois
extensions. The {\em dual} $\mathcal{P}^*$  of $\mathcal{P}$ is the class
 of all duals of extensions in $\mathcal{P}$. 
\end{definition}

In terms of our notational conventions, 
$H$ has $\mathcal{P}$ if and only if $H^*$ has $\mathcal{P}^*$.
\begin{proposition}
\tlabel{dual} Let $\mathcal{P}$ be a property of   finite-dimensional Hopf algebras. If
$\mathcal{P}$  ascends via subnormal series then so does  $\mathcal{P}^*$.  If $\mathcal{P}$
respects field extension/restriction then so does $P^*$.
\end{proposition}
\begin{proof} If $K\to H\to\overline{H}$ is an exact sequence (that is, $K$ is a normal  Hopf
subalgebra of $H$ and $\overline{H}=H/HK^+$), then by dualizing we obtain 
$\overline{H}^*\to H^*\to K^*$ which is also exact, proving the first part.  The second follows
from the fact that for every field $E$ with $E\supset F$, 
 $(H\otimes_F E)^*$ is naturally isomorphic
 to $H^*\otimes_F E$ as a Hopf algebra over $E$.
\end{proof}
\subsection*{Equivalences}In this subsection $H$ will be assumed to be finite-dimensional.
 A basic fact is that $R$ and $A\#H^*$ are Morita equivalent. In the 
following key result from \cite{montgomery-schneider;prime},
$\Phi_1, \Phi_2$  denote the  map $\Phi$ from \ref{fundcorr}  in the extensions $(R,A)$ and
$(A,A\#H^*)$, respectively.

\begin{theorem}
\tlabel{Morita} Let $(R,A)$ be a faithfully flat $H$-Galois extension, with $H$
finite-dimensional. Then the Morita equivalence between $R$ and $A\#H^*$ induces a lattice
isomorphism $f: I\mapsto I^\dagger$ from  $\mathcal{I}(R)$ to
$\mathcal{I}(A\#H^*)$. This map preserves products, and 
$(I^\dagger:H)=\Phi_2\circ\Phi_1((I:H))$.
\noproof 
\end{theorem} Thus given the chain of extensions $R\subset A\subset A\#H^*$ we may freely expand
and contract stable ideals in the natural way.

We wish to develop this correspondence further, into a correspondence between
 extensions.
\begin{theorem}
\tlabel{smashprod} The following diagram commutes. All vertical maps are induced by the Morita
equivalence as in 
\ref{Morita}. The middle horizontal and diagonal maps come from \ref{fundcorr}.
\begin{displaymath}
\begin{diagram}
\node{\mathcal{I}(R)}\arrow{e,b}{(:H)}\arrow{s,<>}
\node{\mathcal{I}_H(R)} \arrow{e,<>} \arrow{s,<>}
\node{\mathcal{I}_{H^*}(A)}\arrow{s,<>}
\node{\mathcal{I}(A)}\arrow{s,<>} \arrow{w,b}{(:H^*)}\\
\node{\mathcal{I}(A\#H^*)}\arrow{e,t}{(:H)}
\node{\mathcal{I}_H(A\#H^*)} \arrow{e,<>}\arrow{ne,<>}
\node{\mathcal{I}_{H^*}(A\#H^*\#H)}
\node{\mathcal{I}(A\#H^*\#H)} \arrow{w,t}{(:H^*)}
\end{diagram}
\end{displaymath}
\end{theorem}
\begin{proof} This all follows from \ref{Morita}. First note  the two central triangles commute
since $(I:H)^\dagger=\Phi_2\circ\Phi_1((I:H))$. This then shows that 
$(I:H)^\dagger=(I:H)(A\#H^*)$ and  so the
left square commutes, as does the right.
\end{proof} It is clear  from \ref{smashprod} that  as far as ideals are 
concerned, the extensions $(R,A)$ and $(A\#H^*,
A\#H^*\#H)$ are equivalent in a strong sense. In order to verify many 
results on ideals in faithfully flat
finite Galois  extensions, it will suffice to consider smash products. 

We also wish to deal with field extension and contraction in a similar way. Let $E$ be an
extension field of $F$, and adopt the notation of \ref{stability}.
\begin{theorem}
\tlabel{fieldext} The following diagram commutes. The maps marked $c$ are contraction, and the
 other vertical maps are induced by field extension and restriction. The  middle horizontal maps
come from \ref{fundcorr}.
\begin{displaymath}
\begin{diagram}
\node{\mathcal{I}(R')}\arrow{s,l}{c}\arrow{e,t}{(:H')}
\node{\mathcal{I}_{H'}(R')}\arrow{s,<>} \arrow{e,<>}
\node{\mathcal{I}_{H^{'*}}(A')}\arrow{s,<>}
\node{\mathcal{I}(A')}\arrow{w,t}{(:H^{'*})}\arrow{s,r}{c}\\
\node{\mathcal{I}(R)}\arrow{e,b}{(:H)}
\node{\mathcal{I}_H(R)}\arrow{e,<>}
\node{\mathcal{I}_{H^*}(A)}
\node{\mathcal{I}(A)}\arrow{w,b}{(:H^*)}\\
\end{diagram}
\end{displaymath}
\end{theorem}
\begin{proof}By the previous result, we may assume that $A=R\#H$ and $A'=R\#H'$. Then  it is
clear from the definition of the action of $H'$ on $R'$ that for each ideal $I$ of $R'$,
$(I:H')\cap R=((I\cap R):H)$,  so the outer squares commute. Also, the middle square commutes
because $(H')^*=(H^*)'$.
\end{proof}

For modules, the situation is similar to the above. The map 
$g:R\text{--Mod}\to R'\text{--Mod}$ which takes $M$ to $M\otimes_F E$ has the property that the
lattice of submodules of $g(M)$ is isomorphic to  the lattice of submodules of $M$. This
property is also  shared by the map induced by Morita equivalence. Both of these maps also
commute with induction and restriction. We summarize the latter fact below.

\begin{theorem}
\tlabel{smashprodmod} The following diagrams commute. The vertical maps are induced by the
Morita equivalence.
\begin{displaymath}
\begin{diagram}
\node{R\textup{--Mod}}\arrow{e,t}{\textup{ind}}\arrow{s,<>}\node{A\textup{--Mod}}\arrow{s,<>}\\
\node{A\#H^*\textup{--Mod}}\arrow{e,t}{\textup{ind}}\node{A\#H^*\#H\textup{--Mod}}
\end{diagram}
\qquad
\begin{diagram}
\node{R\textup{--Mod}}\arrow{s,<>}\node{A\textup{--Mod}}\arrow{w,t} {\textup{res}}\arrow{s,<>}\\
\node{A\#H^*\textup{--Mod}}\node{A\#H^*\#H\textup{--Mod}}\arrow{w,t}{\textup{res}}
\end{diagram}
\end{displaymath} The following diagrams commute. The vertical maps are induced by field
extension.
\begin{displaymath}
\begin{diagram}
\node{R\textup{--Mod}}\arrow{e,t}{\textup{ind}}\arrow{s}\node{A\textup{--Mod}}\arrow{s}\\
\node{R'\textup{--Mod}}\arrow{e,t}{\textup{ind}}\node{A'\textup{--Mod}}
\end{diagram}
\qquad
\begin{diagram}
\node{R\textup{--Mod}}\arrow{s}\node{A\textup{--Mod}}\arrow{w,t} {\textup{res}}\arrow{s}\\
\node{R'\textup{--Mod}}\node{A'\textup{--Mod}}\arrow{w,t}{\textup{res}}
\end{diagram}
\end{displaymath} The following diagrams commute. The vertical maps are induced by field
restriction.
\begin{displaymath}
\begin{diagram}
\node{R\textup{--Mod}}\arrow{e,t}{\textup{ind}}\node{A\textup{--Mod}}\\
\node{R'\textup{--Mod}}\arrow{n}\arrow{e,t}{\textup{ind}}\node{A'\textup{--Mod}}\arrow{n}
\end{diagram}
\qquad
\begin{diagram}
\node{R\textup{--Mod}}\node{A\textup{--Mod}}\arrow{w,t} {\textup{res}}\\
\node{R'\textup{--Mod}}\arrow{n}\node{A'\textup{--Mod}}\arrow{w,t}{\textup{res}}\arrow{n}
\end{diagram}
\end{displaymath}
\end{theorem}
\begin{proof}
The category equivalence between $R$-Mod and $A\#H^*$-Mod is given by $V\mapsto A\otimes_R V$. 
For any Galois extension and any $A\#H^*$-module $M$, there is a natural isomorphism 
$A\otimes_R M^{H^*}\to M$ \cite[8.3.3]{montgomery;hopf-action-rings}. Applying this with 
$M=A\#H^*\otimes_A V$, where $V$ is an $A$-module, yields the second diagram. Now let 
$V$ be an $R$-module and apply the isomorphism with
$A$ replaced by $A\#H^*$, $R$ by $A$ and $M=A\#H^*\#H\otimes_{A\#H^*}(A\otimes_R V)$, yielding
the first diagram. The remaining diagrams are clear from the definition of extended and restricted 
module.
\end{proof}
\subsection*{Prime ideals} An ideal $P$ of $R$ is {\em  $H$-prime} if it is $H$-stable and
whenever $I, J$ are $H$-stable ideals of $R$ with $IJ\subseteq P$ then
$I\subseteq P$ or $J\subseteq P$.

A prime ideal which is $H$-stable is clearly $H$-prime. We write
$H\Spec(R)$ for the poset of all $H$-prime ideals. The equivalence relation $\sim_H$ restricts
to an equivalence relation on  $\Spec(R)$ which we again write as $\sim_H$.

\begin{proposition}
\tlabel{Hprime} Let $(R,A)$ be a faithfully flat $H$-Galois extension.
\begin{enumerate}
\romenumi
\item The maps $\Phi, \Psi$ restrict to poset isomorphisms between
$H\Spec(R)$ and  $H^*\Spec(A)$. We write $I\leftrightarrow J$ in this situation, and say that
$\Phi, \Psi$ respect  prime ideals.
\item The map $Q\mapsto (Q:H)$ is a poset map from $\Spec(R)$ onto
$H\Spec(R)$, and hence induces  a bijection between $\Spec(R)/\sim_H$ and $H\Spec(R)$.
\item The map $P\mapsto P\cap R$ is a poset map from $\Spec(A)$ onto
$H\Spec(R)$, and hence, if $H$ is finite-dimensional, it induces a  bijection between
$\Spec(A)/\sim_{H^*}$ and $H\Spec(R)$.
\end{enumerate}
\noproof
\end{proposition}

For some questions, it suffices to consider the coradical $H_0$ of
$H$. Since $H_0$ is not necessarily a Hopf subalgebra, it is necessary to define the notion of
$C$-stable ideal for an arbitrary subcoalgebra of $H$. Once this is done, as in
\cite{montgomery-schneider;prime}, we can define $C\Spec(R)$ in the obvious way.

\begin{proposition}[\cite{montgomery-schneider;prime}]
\tlabel{H0prime} Let $(R,A)$ be a faithfully flat $H$-Galois extension.
 The map
$I\mapsto (I:H)$ is a poset isomorphism
of $H_0\Spec(R)$ onto $H\Spec(R)$. 
\end{proposition}

\begin{proposition}
\tlabel{specequiv} The diagram in \ref{smashprod} yields the following commutative diagram.
\begin{displaymath}
\begin{diagram}
\node{\Spec(R)}\arrow{e,b}{(:H)}\arrow{s,<>}
\node{H\Spec(R)} \arrow{e,<>} \arrow{s,<>}
\node{H^*\Spec(A)}\arrow{s,<>}
\node{\Spec(A)}\arrow{s,<>} \arrow{w,b}{(:H^*)}\\
\node{\Spec(A\#H^*)}\arrow{e,t}{(:H)}
\node{H\Spec(A\#H^*)} \arrow{e,<>}\arrow{ne,<>}
\node{H^*\Spec(A\#H^*\#H)}
\node{\Spec(A\#H^*\#H)} \arrow{w,t}{(:H^*)}
\end{diagram}
\end{displaymath} The diagram in \ref{fieldext} yields the following commutative diagram.
\begin{displaymath}
\begin{diagram}
\node{\Spec(R')}\arrow{s,l}{c}\arrow{e,t}{(:H')}
\node{H'\Spec(R')}\arrow{s,<>} \arrow{e,<>}
\node{H^{'*}\Spec(A')}\arrow{s,<>}
\node{\Spec(A')}\arrow{w,t}{(:H^{'*})}\arrow{s,r}{c}\\
\node{\Spec(R)}\arrow{e,b}{(:H)}
\node{H\Spec(R)}\arrow{e,<>}
\node{H^*\Spec(A)}
\node{\Spec(A)}\arrow{w,b}{(:H^*)}\\
\end{diagram}
\end{displaymath}
\end{proposition}
\begin{proof}  The bijection $f:I\mapsto I'$ of \ref{Morita} preserves products, and hence 
under $f$ (and its inverse), prime ideals correspond to prime ideals. In the second diagram,
contraction respects prime ideals since the extension $R\subset R'$ is centralizing.
\end{proof}
\section{Krull relations}
\tlabel{sec:lying} For an extension $R\subseteq A$ of commutative rings, there is a well-behaved
relationship between the prime ideals of $A$ and those of $R$, 
given by contraction. This also holds for
finite centralizing extensions. However, in general ring extensions, this correspondence no
longer holds. The standard definition is that $P\in\Spec(A)$
 lies over $Q\in\Spec(R)$ if and only if $Q$ is a minimal prime over $P\cap R$.

In the faithfully flat Galois situation, there are other reasonable definitions.  We can say
that $P$ lies over $Q$ if $P\cap R=(Q:H)$. This is the definition of Montgomery and Schneider,
which we shall use in what follows. Another possibility is: $P$ lies over $Q$ if $P$ is minimal
over $(Q:H)A$.

In the remainder of this section, we assume that $H$ is finite-dimensional  and work with a
fixed faithfully flat $H$-Galois  extension $\mathcal{E}=(R,A)$.  

\begin{definition} Let $P$ be a prime ideal of $A$ and $Q$ a prime ideal of
 $R$. Say that {\em $P$ lies over $Q$} if $(P:H^*)\leftrightarrow (Q:H)$.
\end{definition} 

The above formulation brings out the symmetry in the situation, and we may equally well say
that  $Q$ lies under $P$. More explicitly, $P$ lies over $Q$ if and only if $P\cap R=(Q:H)$, if
and only if
$(Q:H)A=(P:H^*)$.  Note that by \ref{Hprime}, for each  $P\in\Spec(A)$ there is some
$Q\in\Spec(R)$ lying under $P$, and for each $Q\in \Spec(R)$ there is  some $P\in\Spec(A)$ lying
over $Q$.

It follows from \ref{specequiv} that $P\in\Spec{A}$ lies over $Q\in\Spec(R)$ if and only  if
$Q'$ lies over $P$, and this occurs if and only if $P'$ lies over $Q'$. Similarly, 
$P\in\Spec(A')$ lies over $Q\in\Spec(R')$ if and only if $P\cap A$ lies over 
$Q\cap R$. This enables us to streamline considerably the verification of many properties in this
section, since it usually suffices to consider smash products, we can use duality systematically,
and we may extend the field if  necessary.

\subsection*{External Krull relations} Montgomery and Schneider combine the usual  cutting down
and lying over relations into a single stronger property which they call $t$-LO.

\begin{definition}We say that $\mathcal{E}$ has {\em $t$-LO} if and only if
 for every $J\in  H^*\Spec(A)$, there are $P_1,\dots ,P_n\in \Spec(A)$ equivalent to $J$, such
that  $(\cap_i P_i)^t\subseteq J$. We say that $\mathcal{E}$ has {\em $t$-coLO} if for every
$I\in H\Spec(R)$, there are $Q_1,\dots ,Q_m\in \Spec(R)$ equivalent to $I$, such that $(\cap_i
Q_i)^t \subseteq I$. 
\end{definition}  

If we do not wish to specify $t$, then we shall just say that LO or coLO is  satisfied. Note
that the definition in \cite{montgomery-schneider;prime} has the additional requirement that
$m,n\leq \dim H$. We shall not  require such detailed information.
\begin{proposition}
\tlabel{t-LO} Suppose that $\mathcal{E}$ has $t$-LO, let $P\in\Spec(A)$ and let
$Q\in\Spec(R)$. Then 
\begin{enumerate}
\romenumi
\item there are at most $t$ primes of $A$ minimal over a given
$H^*$-prime of $A$
\item if $P$ is minimal over $J\in H^*\Spec(A)$ then $P\cap R=J\cap R$, and 
\item if $P$ is minimal over $(Q:H)A$ then $P$ lies over $Q$.
\end{enumerate} Suppose that $\mathcal{E}$ has $t$-coLO, let $P\in\Spec(A)$ and let
$Q\in\Spec(R)$. Then 
\begin{enumerate}
\romenumi
\item there are at most $t$ primes of $R$ minimal over a given
$H$-prime of $R$
\item if $Q$ is minimal over $I\in H\Spec(A)$ then $(Q:H)=I$, and 
\item if $Q$ is minimal over $P\cap R$ then $P$ lies over $Q$.
\end{enumerate}
\end{proposition}

\begin{proof} We prove only the $t$-coLO case as the other follows by duality. Fix an  $H$-prime
ideal $I$ of $R$  and let $Q'$ be a prime of $R$ minimal over $I$. By $t$-coLO there exist $Q_1,
\dots ,Q_m \in \Spec(R)$ such that $(Q_i:H)=I$ and $(\cap_i Q_i)^t\subseteq I\subseteq Q'$. Since
$Q'$ is prime, some  $Q_i\subseteq Q'$ and by minimality
$Q_i=Q'$. Hence every such $Q'$ is one of the $Q_i$ and $(Q':H)=I$. This establishes (i) and
(ii), and (iii) follows from the way we have defined lying over.
\end{proof}

It is natural to investigate the converse of this last implication.

\begin{proposition}
\tlabel{coINC} The following conditions are equivalent for $Q\in
\Spec(R)$: 
\begin{enumerate}
\romenumi
\item  The prime ideals of $R$ which are $H$-equivalent to $Q$ are mutually incomparable
\item $Q$ is minimal over $(Q:H)$ 
\item For $P\in \Spec(A)$, if $P$ lies over $Q$ then $Q$ is minimal over $P\cap R$ 
\end{enumerate}
\end{proposition}
\begin{proof} Assuming (i), if $(Q:H)\subseteq Q'\subseteq Q$, and $Q'\in \Spec(R)$ is minimal
over  $(Q:H)$, then $Q=Q'$, yielding (ii). If (ii) holds, then if $Q\subseteq Q'$ and
$(Q:H)=(Q':H)=I$, say, then both $Q$ and
$Q'$ are minimal over $I$ and hence are equal. Parts (ii) and (iii) are equivalent by
definition. 
\end{proof}

The dual result is also true.

\begin{proposition}
\tlabel{INC} The following conditions are equivalent for $P\in
\Spec(A)$: 
\begin{enumerate}
\romenumi
\item  The prime ideals of $A$ which are $H^*$-equivalent to $P$ are mutually incomparable
\item $P$ is minimal over $(P:H^*)$
\item  For $Q\in \Spec(R)$, if $P$ lies over $Q$ then $P$ is minimal over $(Q:H)A$
\end{enumerate}
\noproof
\end{proposition}

\begin{definition} Say that $\mathcal{E}$ has {\em coINC} if and only if the conditions in
\ref{coINC} hold for all $Q\in\Spec(R)$. Say that $\mathcal{E}$ has {\em INC} if and only if the
conditions in \ref{INC} hold for all
$P\in\Spec(A)$.
\end{definition} Thus when coINC holds, the map $Q\mapsto (Q:H)$ is strictly increasing, and
when INC holds, the  map $P\mapsto P\cap R$ is strictly increasing. If coINC holds, and $R$ is
prime, then every nonzero ideal of $R$ contains a nonzero $H$-stable ideal. 

Few interesting  results can be achieved in the presence of just one of the Krull  relations. A
useful combination is the property of satisfying  coLO and coINC, which is transitive and holds
over subfields and  algebraic extension fields \cite[Sections 6 and
7]{montgomery-schneider;prime}.  When this property holds,  our definition of lying over is
equivalent to  the standard definition. This is a consequence of the  next Wedderburn-type
result, whose proof is routine (compare \cite[Theorem 16.2] {passman;infinite-crossed}).

\begin{proposition}
$H$ satisfies both coLO and coINC if and only if for every $H$-prime 
$H$-module algebra $R$, there are finitely many primes $Q_1, \dots ,Q_m$ with 
$(Q_i:H)=0$, all the $Q_i$ are incomparable, and their intersection
$N$ is nilpotent. In this  case the $Q_i$ are precisely the minimal primes of $R$ and $N$ is the 
maximum nilpotent ideal of $R$.
\noproof
\end{proposition}

We now come to going up. As noted in \cite{montgomery-schneider;prime}, going down is a 
consequence of LO, and so we shall omit it (but see the next subsection,  where the internal
version is given).

\begin{definition} Say that $\mathcal{E}$ satisfies {\em going up (GU)} if whenever we have
$Q_1, Q_2\in\Spec(R)$ and  $P_1\in\Spec(A)$ such that
$Q_1\subseteq Q_2$ and $P_1$ lies over $Q_1$, then there is
$P_2\in\Spec(A)$  containing $P_1$ and lying over $Q_2$.

Say that $\mathcal{E}$ satisfies {\em co-going up (coGU)} if whenever we have $P_1, P_2\in
\Spec(A)$ and  $Q_1\in \Spec(R)$ such that
$P_1\subseteq P_2$ and $P_1$ lies over $Q_1$, then there is $Q_2\in
\Spec(R)$  containing $Q_1$ and lying under $P_2$.
\end{definition}

In all cases the use of the prefix ``co" is justified --- it is shown in
\cite{montgomery-schneider;prime} that   $t$-LO is dual to $t$-coLO,
 and similarly for the other properties.

\subsection*{Internal Krull relations} We now consider the relationship between $\Spec(R)$ and
$H\Spec(R)$ in more detail, by studying  ``internal" versions of the lying over relations.  The
internal versions of $t$-LO, $t$-coLO, INC and coINC are the same as above. However for going up
and down, there are differences. 

\begin{definition} Say that $\mathcal{E}$ satisfies
 {\em internal coGU} if whenever
 we have $I_1, I_2\in H\Spec(R)$ and $Q_1\in \Spec(R)$ such that
$I_1\subseteq I_2$ and $(Q_1:H)=I_1$ then there is $Q_2\in \Spec(R)$ such that $Q_1\subseteq
Q_2$ and $(Q_2:H)=I_2$. 

Say that $\mathcal{E}$ satisfies  {\em internal GU} if whenever  we have
$Q_1, Q_2\in \Spec(R)$ and $I_1\in \Spec(R)$ such that  $Q_1\subseteq Q_2$ and $(Q_1:H)=I_1$
then there is $I_2\in \Spec(R)$ such that
$I_1\subseteq I_2$ and $(Q_2:H)=I_2$.
\end{definition}

\begin{definition} Say that $\mathcal{E}$  satisfies {\em internal GD} if whenever we have
$Q_1, Q_2\in \Spec(R)$ and  $I_2\in H\Spec(R)$ such that $Q_1\subseteq Q_2$ and $(Q_2:H)=I_2$
then there is $I_1\in H\Spec(R)$  contained in
$I_2$ with $(Q_1:H)=I_1$.

Say that $\mathcal{E}$ satisfies {\em internal coGD} if whenever we have
$I_1, I_2\in H\Spec(R)$ and  $Q_2\in \Spec(R)$ such that $I_1\subseteq I_2$ and $(Q_2:H)=I_2$,
then there is $Q_1\in \Spec(R)$ contained in
$Q_2$ with $(Q_2:H)=I_2$.
\end{definition}

The proof of the following proposition is immediate from the definitions.
\begin{proposition}
\tlabel{int-ext} The following conditions hold for $\mathcal{E}$.
\begin{enumerate}
\romenumi
\item The internal versions of coGU and coGD imply the corresponding external versions
\item The external versions of GU and GD imply the corresponding internal versions
\end{enumerate}
\noproof
\end{proposition}

Note that the dual of internal coGU is not internal GU. In fact it is  easy to see that this
dual condition is formally stronger than GU.

\begin{proposition}
\label{intcoGU}
$\mathcal{E}$ satisfies internal coGU in the following situations:
\begin{enumerate}
\romenumi
\item $\mathcal{E}$ satisfies LO and coGU 
\item $\mathcal{E}$ satisfies coLO, coGU and GU
\item Every prime of $R$ is maximal and $\mathcal{E}$ satisfies either LO or coLO
\end{enumerate}
\end{proposition}

\begin{proof}Suppose that $\mathcal{E}$ satisfies LO and coGU, and let 
$I_1, I_2\in H\Spec(R)$ and $Q_1\in \Spec(R)$ be such that
$I_1\subseteq I_2$ and $(Q_1:H)=I_1$.  Then there is $P\in \Spec(A)$ such that $I_2=P\cap R$. 

By LO, there exist primes $P_1, \dots ,P_n$ of $A$ such that
$P_i\cap R=I_1$ and $(\cap_i P_i)^t\subseteq I_1A\subseteq I_2A\subseteq P$. Since $P$ is prime,
$P_i\subseteq P$ for some
$i$. Now by coGU, there is $Q_2\in \Spec(R)$ so that $Q_2\subseteq Q_1$ and  $P\cap R=(Q_2:H)$.
This proves (i). 

For (ii), again let $I_1, I_2\in H\Spec(R)$ and $Q_1\in \Spec(R)$ be such that  $I_1\subseteq
I_2$ and $(Q_1:H)=I_1$. Write $I_1=P\cap R, I_2=(Q:H)$ for some $Q\in \Spec(R)$ and $P\in
\Spec(A)$. By coLO, there are primes $Q_3,\dots ,Q_m$ of $R$ with $(Q_i:H)=I_1$ and
$(\cap_i Q_i)^t\subseteq I_1\subseteq I_2\subseteq Q$. Thus some
$Q_i\subseteq Q$. Now apply GU to the diagram formed by $Q_i, Q$ and
$P$, to obtain a prime $P'$ of $A$ with $P\subseteq P'$ and $P'\cap R=I_2$.  Applying coGU to
the diagram formed by $I_1, P, P'$ yields
$Q'\in \Spec(R)$ such that  $(Q':H)=P'\cap R=I_2$ and $Q\subseteq Q'$, yielding the result. 
Part (iii) follows from (i) and (ii) because since every prime is maximal, both GU and coGU are
automatic.
\end{proof}

As mentioned in the previous subsection, GD is implied by LO. The dual is  true, and it is
easily seen that the  stronger implication coLO $\Rightarrow$ internal coGD holds.

\subsection*{All Krull relations} The nicest situation is when  all the external Krull relations
hold. Note that in this case all internal relations are also satisfied. It follows in a
straightforward manner from 
\ref{specequiv}, as in
\cite{montgomery-schneider;prime}, that this property is transitive and self-dual, and holds
over all subfields and all algebraic extension fields. This last stipulation is because only for
algebraic extensions can we be guaranteed that $R\subset R'$ has all Krull relations.

If all Krull relations hold,   every strictly ascending (descending) chain 
in any of $H\Spec(R)$,
$\Spec(R)$, $H^*\Spec(A)$  or $\Spec(A)$ yields a strictly ascending 
(descending) chain of the same
 length in each of the other posets. This enables us to compare prime  heights and depths,
classical Krull dimension, etc.

 Group algebras of finite groups, and their duals, satisfy all Krull relations. The proof of this
 rests on work of Lorenz and Passman \cite{lorenz-passman;prime-crossed}. There are no
 Hopf algebras  which are known  not to satisfy all the Krull relations. However, for most
examples, either verifying a Krull relation or else determining that it does not hold is
extremely difficult.
\section{Maximal ideals}
\tlabel{sec:maximal}  As usual, in this section $H$ is a finite-dimensional Hopf algebra and
$\mathcal{E}=(R,A)$ a
 faithfully flat $H$-Galois extension.

We say that $P$ is an {\em $H$-maximal} ideal if
$P$ is a maximal element of  $\mathcal{I}_H(R)$. The poset of all such ideals is denoted by
$H\Max(R)$; it is clearly contained in
$H\Spec(R)$. It is immediate that an $H$-stable ideal which is also maximal is an $H$-maximal
ideal. The quotient of $R$ by an $H$-maximal ideal is an {\em $H$-simple} ring. 

\begin{proposition}
\tlabel{maxcorr} The maps $\Phi$ and $\Psi$ respect maximal ideals; that is, they restrict to
poset isomorphisms $H\Max R\leftrightarrow H^*\Max A$. 
\end{proposition}
\begin{proof} The set of all $H$-stable ideals containing a given $H$-stable
 ideal is inductive and so by Zorn's lemma, every $H$-stable ideal is contained in an
$H$-maximal ideal. Thus an $H$-maximal ideal is precisely a maximal element of $H\Spec(R)$.
Since we know that $\Phi$ and $\Psi$ restrict to maps between $H\Spec(R)$ and $H^*\Spec(A)$, and
the set of maximal elements is a poset invariant, the result follows.
\end{proof}

We  wish to see how the lying over relations are compatible with maximal ideals.

\begin{proposition}
\tlabel{maxLO} The following conditions hold for $\mathcal{E}$.
\begin{enumerate}
\romenumi
\item Let $I$ be an $H$-maximal ideal of $R$. Then there is a maximal ideal $M$ of $A$  with
$M\cap R=I$.
\item  Let $I$ be an $H$-maximal ideal of $R$. Then there is a maximal ideal $M$ of $R$ with 
$(M:H)=I$. 
\end{enumerate}
\end{proposition}

\begin{proof}Let $M$ be a maximal ideal of $A$ (respectively, of
$R$) containing $IA$ (respectively, $I$).
\end{proof}

Consider the property that in \ref{maxLO}(i), {\em every} prime ideal of $A$ with $P\cap R=I$ is
maximal. This property can be considered as a weakened form of  INC, since INC  implies it by
\ref{maxLO}. Similarly, if coINC holds, then every prime ideal of $R$ with core $I$ is maximal.
Thus coINC implies that every prime $H$-simple ring is in fact simple. 

If the  weakened forms of INC and coINC hold then there is no confusion when speaking of the
equivalence class of a maximal ideal, as the equivalence classes in $\Spec$ and $\Max$ coincide.
A priori,  however, there is no reason to expect these two equivalence classes to be equal.

Now by \ref{maxLO} there is a map $H\Max(R)\to \Max(R)/\sim_H$ which is $1-1$. Similarly there
is a map $H\Max(R)\to\Max(A)/\sim_{H^*}$ which is $1-1$. In contrast to the situation with
$\Spec$, it is not clear that these maps are onto. The following definitions address this
question.

\begin{definition} We say that $\mathcal{E}$ has {\em coMAX} if  whenever
$M$ is a maximal ideal of $R$ then $(M:H)$ is an $H$-maximal ideal of
$R$. Dually, $\mathcal{E}$ has {\em MAX} if whenever $M$ is a maximal ideal of $A$ then $M\cap
R$ is an $H$-maximal ideal of $R$. 
\end{definition}

It is readily seen, using \ref{maxcorr}, that MAX and coMAX are indeed dual to each other. Now
suppose that $K$ is a normal Hopf subalgebra of $H$ with quotient
$\overline{H}$, and consider the chain $R\subset B\subset A$ as in Section~\ref{sec:basics}. 
Suppose that each intermediate extension has both MAX and coMAX. If $M$ is maximal in $A$, then
$M\cap B$ is $H$-maximal in $B$. By \ref{maxLO}, there is a maximal ideal
 $M'$ of $B$ with $(M':\overline{H})=M\cap B$. It follows that $M\cap R=((M'\cap R):H)$. This
latter ideal is $H$-maximal by MAX and coMAX applied to $(R,B)$. Thus by \ref{dual}, the property
of satisfying MAX and coMAX is transitive.

The proof of the following proposition is immediate from the definitions.

\begin{proposition}
\tlabel{coGU-coMAX}  If $\mathcal{E}$ satisfies internal coGU then
$\mathcal{E}$ satisfies coMAX.  If $\mathcal{E^*}$ satisfies internal coGU 
then $\mathcal{E}$ satisfies MAX.  \noproof
\end{proposition}

\begin{theorem}
\tlabel{max} Suppose that $\mathcal{E}$ satisfies coINC and coMAX. Then following conditions are
equivalent for $Q\in \Spec(R)$.
\begin{enumerate}
\romenumi
\item $Q$ is maximal 
\item $(Q:H)$ is $H$-maximal
\item $Q$ is minimal over an $H$-maximal ideal of $R$
\item $(Q:H)A$ is $H^*$-maximal
\end{enumerate}

Dually, suppose that $\mathcal{E}$ satisfies INC and MAX. Then the following conditions are
equivalent for $P\in\Spec(A)$.
\begin{enumerate}
\romenumi
\item $P$ is maximal
\item $(P:H^*)$ is $H^*$-maximal
\item $P$ is minimal over an $H^*$-maximal ideal of $A$
\item $P\cap R$ is $H$-maximal
\end{enumerate}
\end{theorem}

\begin{proof} We prove only the equivalence of the first four conditions. Clearly (ii) and (iv)
are equivalent by \ref{maxcorr}. Also (ii) implies (iii) by \ref{coINC}, and given (iii), the
$H$-maximal ideal in question must be
$(Q:H)$, so (ii) and (iii) are equivalent. Now (ii) implies (i) in the presence of coINC and the
implication (i)$\Rightarrow$(ii) is precisely coMAX. 
\end{proof}

\begin{corollary}
\tlabel{maxLOmax}  Consider the property of $H$: if $P\in \Spec(A)$ lies over $Q\in
 \Spec(R)$, then $P$ is maximal if and only if $Q$ is maximal.

This property is self-dual, transitive and field-independent, and it holds
 if $H$ satisfies MAX, coMAX, INC and coINC.
\noproof
\end{corollary} The property in \ref{maxLOmax} can be considered as a minimal requirement for a
decent theory of lying over. Recall that as far as we know, every  finite-dimensional $H$
satisfies the hypotheses of \ref{maxLOmax}.

\section{Modules}
\tlabel{sec:modules} Let $\mathcal{E}=(R,A)$ be a ring extension, $V$ an $R$-module, and $W$ an
 $A$-module. Then we
 denote the restricted $R$-module  $W_{\mid R}$ by $W^\d$ and the induced
 $A$-module $A\otimes_R V$ by $V^\u$. 

\begin{proposition}
\tlabel{ann} Suppose that $\mathcal{E}$ is faithfully flat $H$-Galois, and  let $V$ be an
$R$-module and $W$ an $A$-module. Then the following conditions
 hold.
\begin{enumerate}
\romenumi
\item $\ann W^\d = (\ann W)\cap R$
\item $\ann V^{\u\d} = (\ann V:H)$
\item If $H$ is finite-dimensional, then $\ann V^{\u} = (\ann V:H)A$
\item If $H$ is finite-dimensional, then $\ann W^{\d\u}=(\ann W:H^*)$
\item If $V$ is simple and $X$ is a simple image of $V^\u$, then $\ann X$ lies over $\ann V$ 
\item If $W$ is simple and $X$ is a simple image of $W^\d$, then $\ann W$ lies over $\ann X$
\end{enumerate}
\end{proposition}

\begin{proof} Part (i) is trivial. For (ii), let $Q=\ann V, P=\ann V^{\u}$.
 Then $P\cap R\subseteq Q$ since $V^{\u\d}$ has an $R$-submodule isomorphic
 to $V$. Since $P\cap R$ is the intersection of an ideal of  $A$ with $R$,  it is $H$-stable and
so $P\cap R\subseteq (Q:H)$. But clearly $(Q:H)AV= A(Q:H)V=0$ so that $(Q:H)\subseteq P\cap R$.
This yields (ii).  Parts (iii) and (iv) follow from the observation that $\ann V^\u$, as the
annihilator of an $H^*$-stable module, is $H^*$-stable.

We now prove (v). Since $X$ is simple, $X=V^\u/M$ for some maximal
 $A$-submodule  $M$ of $V^\u$. Let $P=\ann X, Q=\ann V$. Then
 $(Q:H)A\subseteq P$ by  (ii).  Also $(P\cap R)V$ is an $R$-submodule
 of $V$. Since $(P\cap R)V^\u \subseteq  M\neq V^\u$ and $P\cap R$
 commutes with $A$ it follows that $(P\cap R)V\neq V$  and so $(P\cap
 R)V=0$, since $V$ is a simple $R$-module. Thus $P\cap R\subseteq Q$
 and since  $P\cap R$ is $H$-stable, (v) follows.

Part (vi) follows in a similar manner. If $W$ is simple, let $P=\ann W$,  write $X=W^\d/M$ for
some maximal $R$-submodule $M$ of $W^\d$, and  let $Q=\ann V$. Clearly $P\cap R\subseteq Q$ and
so $P\cap R\subseteq (Q:H)$.  Conversely $(Q:H)W$ is an $A$-submodule of $W$ and so since 
$(Q:H)W\subseteq M\neq W$, we have $(Q:H)W=0$. Thus  the reverse  containment is shown and we
obtain $(Q:H)=P\cap R$. 
\end{proof}

\begin{definition} Say that 
$\mathcal{E}=(R,A)$ has the {\em finite induction
 property} if  whenever $V$ is an $R$-module of finite length, then 
$V^\u$ is an $A$-module of finite length.
\end{definition}

It is easy to see that $H$ has the finite induction property if and only if $H^*$ has the dual
{\em finite restriction property}, that the restriction of  each finite length $A$-module to $R$
is of finite length.

Clearly, both the finite induction and finite restriction  properties are transitive.  It
follows from \ref{smashprodmod} that they are field-independent.
 They are both satisfied by group algebras,  and the finite induction 
property by restricted enveloping algebras 
\cite[Lemma 23]{chin;prime-differential} (we reproduce a proof of a generalization 
of this last result in section~\ref{sec:examples}). In fact these Hopf 
algebras have the stronger property that if $V$ has finite length, then so 
does $V^{\u\d}$.

\begin{definition} Say that $\mathcal{E}=(R,A)$ has the {\em semisimple induction property} if
 whenever $V$ is semisimple of finite  length then so is
 $V^\u$. Dually, $\mathcal{E}$ has the {\em semisimple restriction
 property} if  whenever $W$ is semisimple of finite length then so is $W^\d$.
\end{definition}

Again, the semisimple induction and restriction properties are transitive and 
field-independent. The latter property is always satisfied by group algebras, as is the
former in characteristic coprime to the group order. The semisimple  restriction property is in
fact satisfied by finite normalizing
 extensions (see, for example, \cite{passman;prime-normalizing}). In addition,  such extensions
satisfy the property that $V^{\u\d}$ is semisimple of finite length
 whenever $V$ is.

\section{Primitive ideals}
\tlabel{sec:primitive}  In this section, $H$ is finite-dimensional and $\mathcal{E}=(R,A)$ a 
fixed faithfully flat Galois $H$-extension.

It is natural to define the set $H\Prim(R)$ of {\em $H$-primitive}
 ideals. Unlike the case for prime and maximal ideals, there does not
 seem  to be any definition which is internal to $R$. In the case
 where $A=R\#H$ is a smash product, the obvious definition is that $Q$
 is $H$-primitive if and only if $Q$ is the annihilator of an
 $H$-stable module.  The only reasonable definition of $H$-stable
 module for $R$ is a module for $R\#H$. Thus $Q$ should be
 $H$-primitive if and only if $Q=P\cap R$ for some $P\in \Prim(A)$. 

However, in more general extensions, this notion of $H$-stable module
 does not make sense, since there need be no action of $H$ on
 $R$. Any definition of $H$-primitive should  satisfy the following
 requirements. First, we should have  $H\Max  R\subseteq
 H\Prim(R)\subseteq H\Spec(R)$. Second, an $H$-stable ideal which is
 primitive should be $H$-primitive. Third, the bijections between
 $H\Spec(R)$ and $H^*\Spec(A)$ should yield bijections between
 $H\Prim(R)$ and $H^*\Prim(A)$.

The definition we choose yields the desired results fairly
 quickly, and in order to obtain a primitive analogue of
 \ref{specequiv}, it is the obvious one.
 However, some problems remain. In order to obtain an analogue of
 \ref{H0prime}, it would be necessary to define
 $C$-primitive ideals for an arbitrary subcoalgebra of $H$. No obvious
 candidate presents itself.

\begin{definition}  An ideal $I$ of $R$ is {\em $H$-primitive} if  $I=P\cap R$ for some
 primitive ideal $P$ of $A$. 
\end{definition}

Clearly every $H$-primitive ideal is $H$-prime. Let $I$ be an
 $H$-maximal ideal of $R$. Then by \ref{maxLO}, $I=M\cap R$ for some maximal ideal $M$
 of $A$, so that $I$ is $H$-primitive.  Thus the first  requirement of
 our definition is satisfied. The second and third follow from the
 next result, which shows that another candidate for the definition of
 $H$-primitive is equivalent to the one given. Thus, the intersections
 with $R$ of primitive ideals of $A$ coincide with the $H$-cores of
 primitive ideals of $R$.

\begin{proposition}
\tlabel{Hprim}  The following conditions hold for $\mathcal{E}$.
\begin{enumerate}
\romenumi
\item The map $Q\mapsto (Q:H)$ is a poset epimorphism from $\Prim(R)$ onto  $H\Prim(R)$. This
induces  a bijection between $\Prim(R)/\sim_H$ and $H\Prim(R)$.
\item  The map $P\mapsto P\cap R$ is a poset epimorphism from
 $\Prim(A)$ onto $H\Prim(R)$. This  induces a  bijection between
 $\Prim(A)/\sim_{H^*}$ and $H\Prim(R)$.
\item The poset isomorphism $H\Spec(R)\leftrightarrow H^*\Spec(A)$
 respects primitive ideals; that is, it restricts to a  poset
 isomorphism $H\Prim(R)\leftrightarrow H^*\Prim(A)$. 
\end{enumerate}
\end{proposition}

\begin{proof}Part (ii) is immediate from the way we have defined
 $H$-primitive  ideals. We prove (i). We first show that the map is well defined. 
Let $Q\in \Prim(R)$. Then $Q=\ann V$ for some simple module $V$. Then 
$V^\u$ is finitely generated, hence has a maximal submodule and therefore a simple 
quotient.  Applying \ref{ann}(v) we obtain $P\in \Prim(A)$ such that 
$P\cap R=(Q:H)$. We now show that the map is onto. Given $I\in H\Prim(R)$, write
 $I=P\cap R$ where $P\in \Prim(A)$. Then $P=\ann W$ for some simple 
$A$-module $W$. Since $A$ is a finitely generated $R$-module, $W^\d$ is finitely 
generated. Similarly to the above, applying \ref{ann}(vi) we obtain 
$Q\in \Prim(R)$ with $(Q:H)=I$.  Thus the map is onto, proving (i). Part (iii)
 is now immediate.
\end{proof}

In particular, this last result shows that the analogue for primitive
 ideals of coMAX is always satisfied (the analogue of MAX is satisfied by
 definition).

Since the ideal equivalence of \ref{Morita} is a composition of two 
equivalences of  (bi)module categories and hence preserves
exact sequences, it follows that under this  map, primitive ideals
 correspond to primitive ideals. This, and  the fact that primitive ideals 
behave well with respect to lying over in centralizing extensions, shows that  
the analogue for Prim of \ref{specequiv} holds, if we restrict all maps from Spec
 to Prim. 

It is clear that each of the Krull relations of section~\ref{sec:lying}  can be formulated for
$\Prim$, so that we can speak of, for  example, $t$-LO for primitive ideals. Obviously INC and
coINC  for $\Spec$ imply the corresponding relations for $\Prim$, but the  situation with regard
to the other relations is not immediately apparent. However, a construction of Passman, the
``primitivity machine" 
\cite{passman;prime-normalizing}, can be adapted to our situation, as we shall now explain. 

For every ring $R$, there is an overring $S=\widehat{R}=R\langle X\rangle$
 which is obtained  from $R$ by taking (noncommutative) formal power series and polynomials. The 
variables adjoined all commute with $R$. 
Key properties are: $\widehat{I}\cap R=I, \widehat{R}/\widehat{I}\cong
 \widehat{R/I}, \widehat{A}\widehat{B}\subseteq \widehat{AB}, 
\widehat{\cap_i  A_i}=\cap_i \widehat{A_i}$. Here $I$ is an ideal of $R$
 and $A$ a subset of $R$. If $Q$ is a prime ideal of $R$ then $\widehat{Q}$ is a primitive ideal
of 
$\widehat{R}$. Furthermore, the map $\widehat{}$ is injective on ideals and if $I$ is an annihilator
ideal of $\widehat{R}$, then $I\cap R$ is a prime ideal of $R$.

If $H$ acts on $R$, then the action extends naturally to an action on 
$S$, by letting $H$ act trivially on the adjoined variables in $X$. Furthermore it is clear that
if $A=R\#H$ then $\widehat{A}=\widehat{R}\#H$. We summarize the main facts on spectra in the 
following result, whose proof follows directly from \ref{Hprime}, \ref{Hprim} and the remarks 
above.

\begin{theorem}
\tlabel{SpecPrim} The following diagram commutes. The maps $\iota$ are the natural inclusions,
 and the vertical maps come from \ref{Hprime} and \ref{Hprim}.

\begin{displaymath}
\begin{diagram}
\node{\Prim(A)}\arrow{e,t}{\iota}\arrow{s,l}{(:H^*)}
\node{\Spec(A)}\arrow{e,t}{\widehat{}}\arrow{s,l}{(:H^*)}
\node{\Prim(\widehat{A})}\arrow{s,r}{(:H^*)}\\
\node{H^*\Prim(A)}\arrow{e,t}{\iota}\arrow{s,<>}
\node{H^*\Spec(A)}\arrow{e,t}{\widehat{}}\arrow{s,<>}
\node{H^*\Prim(\widehat{A})}\arrow{s,<>}\\
\node{H\Prim(R)}\arrow{e,t}{\iota}
\node{H\Spec(R)}\arrow{e,t}{\widehat{}}
\node{H\Prim(\widehat{R})}\\
\node{\Prim(R)}\arrow{e,t}{\iota}\arrow{n,l}{(:H)}
\node{\Spec(R)}\arrow{e,t}{\widehat{}}\arrow{n,l}{(:H)}
\node{\Prim(\widehat{R})}\arrow{n,r}{(:H)}\\
\end{diagram}
\end{displaymath}
\noproof
\end{theorem}

It is a consequence of \ref{SpecPrim} that each Krull relation for  Prim implies the analogous
one for Spec.  We give one example here. Suppose that $H$ satisfies INC for Prim. Let $P_1,
P_2\in \Spec(A)$ with $P_1\subset P_2$. Then $\widehat{P_1}\subset
\widehat{P_2}$,  and so since INC holds for Prim, we have $\widehat{P_1}\cap \widehat{R}\subset
\widehat{P_2}\cap\widehat{R}$. Thus we must have $P_1\cap R\subset P_2\cap R$ and so INC holds
for Spec. See the discussion of cocommutative Hopf algebras in  section~\ref{sec:examples}, or
\cite{passman;prime-normalizing} for more details. 

We now investigate the converse implications.
  
\begin{proposition}
\tlabel{t-LO-FIP} If $\mathcal{E}$ satisfies the finite induction property (respectively,  the
finite restriction property), then it satisfies LO (respectively coLO) for primitive ideals.
\end{proposition}

\begin{proof} We prove only the second assertion as the first then follows by duality. Let $P$
be a primitive ideal of $A$, with $P=\ann W$. By \ref{ann}, 
$\ann W^\d=P\cap R$, an $H$-primitive ideal of $R$. Let $Q_1, \dots ,Q_m$ be the annihilators of 
the finitely many composition factors of $W^\d$, where
 $Q_i, 1\leq i\leq s\leq m$, are the distinct ones. Then  whenever $t$ is such that $st\geq m$,
$(\cap_i Q_i)^t\subseteq \ann W^\d=P\cap R$ as required.
\end{proof}

The proof of the following proposition is essentially the same as  the argument of \ref{t-LO}.
\begin{proposition}
\tlabel{minimalprim}  The following conditions hold for $\mathcal{E}$.
\begin{enumerate}
\romenumi
\item Suppose that $\mathcal{E}$ has LO for primitive ideals. If $Q\in
 \Prim(R)$ and $P\in \Spec(A)$  is minimal over $(Q:H)A$ then $P$ is
 primitive.
\item Suppose that $\mathcal{E}$ has coLO for primitive ideals. If $P\in
 \Prim(A)$ and $Q\in  \Spec(R)$ is minimal over $P\cap R$ then $Q$ is
 primitive. 
\noproof
\end{enumerate}
\end{proposition}

\begin{theorem}
\tlabel{primlem}  If $\mathcal{E}$ has  coINC and coLO  for primitive ideals, then the following
conditions are equivalent for $Q\in\Spec(R)$.
\begin{enumerate}
\romenumi
\item $Q$ is primitive 
\item $(Q:H)$ is $H$-primitive
\item $Q$ is a minimal prime over an $H$-primitive ideal of $R$
\item $(Q:H)A$ is $H^*$-primitive
\end{enumerate}

 If $\mathcal{E}$ has INC and LO   for primitive ideals, then the following conditions are
equivalent for 
$P\in \Spec(A)$.
\begin{enumerate}
\romenumi
\item $P$ is primitive 
\item $(P:H^*)$ is $H$-primitive
\item $P$ is a minimal prime over an $H^*$-primitive ideal of $A$
\item $P\cap R$ is $H$-primitive
\end{enumerate}
\end{theorem}

\begin{proof} We prove only the equivalence of the first four conditions, the others  following
by duality. The implication (i)$\Rightarrow$ (ii) is \ref{Hprim}(i), and
 (ii)$\Rightarrow$ (iii) follows from  coINC. Also (ii) and (iv) are
 equivalent by \ref{Hprim}(iii). Now (iii) implies (i) by
 \ref{minimalprim} and \ref{Hprim}(ii). 
\end{proof}
\begin{theorem}
\tlabel{primLOprim} 
 Consider the following property of $H$: if $P\in \Spec(A)$ lies over
 $Q\in\Spec(R)$, then $P$ is primitive if and only if $Q$ is primitive (``lying over respects 
primitivity"). 

This property is self-dual, transitive, and field-independent, and it holds if (i) $H$ satisfies
INC, coINC, LO and  coLO for primitive ideals or (ii) $H$ and $H^*$ are pointed. Conversely,
if lying over respects primitivity and $H$ has all Krull relations for $\Spec$, then 
$H$ has all Krull relations for $\Prim$.
\begin{proof}
The property is self-dual, transitive and field-independent by \ref{specequiv} and the remarks
after \ref{Hprim}.
It holds in case (i) by \ref{primlem}. In case (ii), we use the fact (see
Section~\ref{sec:examples}) that for every pointed $H$,
the primes minimal over a given $H$-prime ideal of $R$ are all conjugate under $G(H)$, so that if
one of them is primitive,  then they all are. By \ref{Hprim} there is one such primitive ideal.
The case where $H^*$ is pointed follows by duality. If $H$ has all Krull relations for Spec and
lying over respects primitivity, then if we start with primitive ideals, all prime ideals produced
by LO, GU and their duals are primitive, as required. We already know that INC and coINC are
inherited by Prim from Spec. 
\end{proof}
\end{theorem}

 As with \ref{maxLOmax}, the property discussed in \ref{primLOprim} can be considered a minimum
requirement of a decent theory of lying over. It holds for finite normalizing extensions
 \cite{mcconnell-robson;noncommutative-noetherian}, and the extensions studied by E.~Letzter  in
\cite{letzter;primitive-noetherian}, namely ring extensions $R\subseteq S$ where $R$ and $S$ are
noetherian and $S$ has finite GK-dimension over $R$ on the left and the right. 

Suppose that all Krull relations for Prim are
satisfied (so that also all Krull relations for Spec hold, and lying over respects primitivity).
This  enables us to prove many results relating Spec and Prim. For example, $A$ satisfies the
property that every primitive ideal is maximal if and  only if $R$ does. The result was first
established for group algebras in \cite[Theorem 1.7]{lorenz;primitive-crossed}. 

\section{Extensions with a total integral}
\tlabel{sec:totint}The faithfully flat Galois property does not always  hold for $H$-extensions
which arise in practice, and can be difficult to  verify. As was the case in
\cite{montgomery-schneider;prime}, 
 some of our results will extend to a larger class of $H$-extensions,
 namely those with a total integral. We shall not go into great detail, but  confine our
discussion to the  properties of such extensions which enable a reasonably systematic
 translation of results on faithfully flat Galois extensions to the more  general context.

The key facts about such an $H$-extension $(R,A)$ are as follows. See
\cite{montgomery-schneider;prime} for more details.
 In (iv), the lying over relation is the  standard one.
\begin{proposition}
\tlabel{totint} Suppose that $H$ is finite-dimensional, and that $(R,A)$ is an $H$-extension
 with a total integral. Then the following statements hold.
\begin{enumerate}
\romenumi
\item There is an idempotent $e\in S=A\#H^*$ such that $eSe=Re\cong R$.
\item This induces a lattice 
map $\mathcal{I}(S)\to \mathcal{I}(R)$ which in turn  induces poset
epimorphisms $\Spec(S)\to\Spec(R)$ and $\Prim(S)\to\Prim(R)$.
\item These latter maps are such that $ePe=R$ if and only if $e\in P$, and they yield poset
isomorphisms $\Spec_e(S)\equiv\{P\in\Spec(S)|e\not\in P\}\to\Spec(R)$ and
$\Prim_e(S)\to\Prim(R)$. 
\item If $H$ has coLO and coINC, then $P$ lies over $eQe$ if and only if $Q$ lies over $P$.
\end{enumerate}
\noproof
\end{proposition}

If $H$ is finite-dimensional and cosemisimple, then every $H$-extension has a  total integral
(see \cite[section 4.3]{montgomery;hopf-action-rings}). On the other hand, if $H$ is connected
then an extension with a total integral is necessarily a crossed product extension
\cite{bell;comodule}.

It follows in a straightforward manner from \ref{totint} (see
\cite[Section 5]{montgomery-schneider;prime}) that if $H$ has all Krull relations for
$\Spec$, then $(R,A)$ inherits the analogues of LO, coLO, INC and coINC. Similar results can be
 obtained for $\Prim$ and the details are left to the reader. However, it is known that not all 
Krull relations hold for such extensions; an example of Montgomery and Small
\cite{montgomery-small;integrality-prime} shows that the analogue of coGU fails even for the
group algebra of a group of order $2$. 

We give one simple example of how \ref{totint} can be applied. Let $R$ be the universal
enveloping algebra of a finite-dimensional nilpotent Lie algebra over a field $F$ of
characteristic zero, and let $G$ be a finite group of automorphisms of $R$. Then it is well
known that  every primitive ideal of $R$ is maximal.  Hence the same is true for $R\#G$. We
claim that the same is also true of $R^G$. The extension $R^G\subseteq R$ has a total integral
since $FG$ is semisimple by Maschke's theorem. Let $P$ be a primitive ideal of $R^G$. Then
$P=eP'e$ for some primitive ideal of $R\#G$. But then $P'$ is maximal and thus so is $P$.

\section{Strongly semiprimitive Hopf algebras}
\tlabel{sec:radicals} We begin by defining the various obvious notions of radical. As usual, $H$
is a Hopf algebra which is  not assumed to be finite-dimensional. We denote
 the  prime radical by $N(R)$ and the Jacobson radical by $J(R)$.

\begin{definition} Let $(R,A)$ be a faithfully flat $H$-extension.  The {\em $H$-prime radical}
$N_H(R)$ of $R$ is the intersection of all 
$H$-prime ideals of $R$.  The {\em $H$-Jacobson radical} $J_H(R)$ of $R$ is the intersection of
all $H$-primitive ideals of $R$.
\end{definition}

\begin{proposition} 
\tlabel{Hrad} The following conditions hold for a faithfully flat
 $H$-Galois extension $(R,A)$.
\begin{enumerate}
\romenumi
\item $J_H(R)=J(A)\cap R=(J(R):H)$
\item $N_H(R)=N(A)\cap R=(N(R):H)$
\end{enumerate}
\end{proposition}

\begin{proof} In (i), the first equality follows immediately from the definition of 
$H$-primitive, and the second from  \ref{Hprim}. In (ii), both parts follow from  \ref{Hprime}.
\end{proof}
 
Because $J_H(R)=(J(R):H)$, the fact that $J_H$ really is a radical in the usual sense   (for,
say, the class of all $H$-module algebras) follows readily from the fact that $J$ is a radical. 
A complete treatment would characterize $J_H$ in terms of  maximal $H$-stable left (or right)
ideals and  left (or right) quasiregularity, and $N_H$ in terms of an $H$-stable analogue of the
transfinite Baer process for
$N(R)$. We shall not pursue this here as our main interest is in algebras with zero radical.

 The following proposition  summarizes, in our notation, several  results in 
\cite[Section 8]{montgomery-schneider;prime}. We say that two conditions are  equivalent on a
class of extensions if whenever one condition holds for  all extensions in the class, then the
other hold for all extensions.  It is not necessarily the case that both conditions are
equivalent for a  fixed extension in the class.

\begin{proposition}
\label{scss}  The following conditions are equivalent on the class of all faithfully flat
$H$-Galois extensions $(R,A)$.
\begin{enumerate}
\romenumi
\item If $R$ is $H$-semiprime then it is semiprime
\item $N(R)$ is $H$-stable
\item $N_H(R)=N(R)$
\item $N(R)\subseteq N(A)$
\item If $A$ is semiprime then $R$ is semiprime
\end{enumerate} Suppose now that $H$ is finite-dimensional and  that these conditions do hold
for all faithfully flat $H$-Galois extensions.
 Then (iv) and (v) in fact hold for all $H$-extensions.
\noproof
\end{proposition}

In the situation of the last paragraph of 
 \ref{scss}, then  we say that $H$ is {\em strongly cosemiprime}. When 
$H$ is finite-dimensional, the dual
 property, strong semiprimeness, is also of interest. $H$ is strongly  semiprime if and only if
whenever $R$ is an $H$-semiprime $H$-module algebra,  then $R\#H$ is semiprime. This is of
course equivalent to many other  conditions dual to those in \ref{scss}. Both properties are
transitive and field-independent.

We move on to primitive ideals, with the aim of obtaining analogous results.

\begin{proposition} The following conditions are equivalent on the class of all faithfully flat 
$H$-Galois extensions $(R,A)$.
\begin{enumerate}
\romenumi
\item If $R$ is $H$-semiprimitive then it is semiprimitive
\item $J(R)$ is $H$-stable
\item $J_H(R)=J(R)$
\item $J(R)\subseteq J(A)$
\item If $A$ is semiprimitive then $R$ is semiprimitive
\end{enumerate} Suppose that $H$ is finite-dimensional and that all conditions hold for all
faithfully flat $H$-Galois extensions. Then
 (iv) and (v) in fact hold for all $H$-extensions.
\end{proposition}
\begin{proof} Clearly (v) and (i) are equivalent by \ref{Hprim}.  Also (ii), (iii) and (iv) are
equivalent for a given extension by \ref{Hrad}.  The remaining equivalences are obtained by
reducing  to the faithfully flat $H$-Galois extension $(R/J), A/AJ)$, where 
$J=J_H(R)$.

Suppose that all the conditions hold for every faithfully flat $H$-Galois
 extension. It suffices to prove that (v) holds for all extensions. Applying (v) in the case
$H^*\subseteq H^*\#H$, we see that $H$ is cosemisimple and so every
$H$-extension has a total integral. If now $A$ is semiprimitive then so is
 $S=A\#H^*$. Thus  there are primitive ideals $P_i$ of $S$ with $\bigcap_i P_i=0$. It
 follows from \ref{totint} that $eP_ie$ is either $R$ or a primitive ideal
 of $R$. Thus deleting those which equal $R$ we obtain a set of primitive  ideals of $R$ with
intersection $e0e=0$, and so $R$ is semiprimitive. 
\end{proof} A Hopf algebra $H$ satisfying (i)-(v) above for all  faithfully flat Galois
$H$-extensions shall be called {\em strongly  cosemiprimitive}.  The characterization of the
dual property, strong semiprimitivity, is left  to the reader. Both properties are transitive
and field-independent.

 Group algebras are certainly strongly cosemiprimitive  (since the Jacobson radical is a
characteristic ideal), and whenever they are semisimple (that is, in characteristic coprime to
the group order),  they are strongly semiprimitive \cite{villamayor;semisimplicity}.

We now relate strong semiprimitivity to other important properties.
\begin{theorem}
\tlabel{ssprim} Let $H$ be a finite-dimensional Hopf algebra and let 
$\mathcal{E}=(R,A)$ be a faithfully flat $H$-Galois extension.
\begin{enumerate}
\romenumi
\item If $\mathcal{E}$ satisfies the semisimple induction
 property then it satisfies 1-LO for primitive ideals.  
\item If $H$ satisfies 1-LO for primitive ideals
 then $H$ is strongly semiprimitive.
\item If $H$ is strongly (co)semiprimitive then $H$ is strongly (co)semiprime.
\end{enumerate}
\end{theorem}
\begin{proof} Let $J$ be an $H^*$-primitive ideal of $A$. By \ref{Hprim}, $J=(Q:H)A$ where 
$Q$ is a primitive ideal of $R$. Let $V$ be a simple $R$-module with  annihilator  $Q$. Then
$V^\u$ is semisimple of finite length and by 
\ref{ann}(iii) its annihilator is $J$. Let $P_1, \dots,  P_n$ be the annihilators of the
composition factors of $V^\u$. By \ref{ann}(v),  each  $P_i$ lies over $Q$. Also $\cap_i P_i$
annihilates $V^\u$ and so this  intersection is contained in $J$. This proves (i).  Now (ii)
follows by taking $J=0$ and using the conclusion of (i).  For (iii) it suffices by duality to
prove only that strongly semiprimitive  implies strongly semiprime.  For this, we use the
primitivity machine. If $A=R\# H$ is semiprime then $\widehat{A}$ is semiprimitive, so that
$\widehat{R}$ is semiprimitive and hence $R$ is semiprime, yielding the desired  result.
\end{proof}

\section{Special classes of Hopf algebras}
\tlabel{sec:examples}
\subsection*{Pointed Hopf algebras}
\tlabel{pointed} Let $H$ be a finite-dimensional pointed Hopf algebra, with coradical filtration of length $t$ and $\dim H_0=d$.  Then it is known
\cite[4.10]{montgomery-schneider;prime} that $H$ satisfies coINC, $t$-coLO, GU
 and coGU.  It follows from \ref{intcoGU} that $H$ and $H^*$ satisfy internal
 coGU, and hence both MAX and coMAX by \ref{coGU-coMAX}. 

We now show that $H$ satisfies the finite induction property.
 In fact this follows from \cite[1.5, 2.1]{schneider;representation-hopf}, and I 
thank Schneider for pointing this out to me.
Those results state that for every simple $R$-module $V$, the module $V^{\u\d}$
 has a filtration $0=V_0\subset V_1\cdots \subset V_n=V$ with finitely many
 terms. Furthermore each quotient is a finite direct sum of modules $V_g$, 
where $g\in G(H)$.  These $V_g$ are all simple $R$-modules and so $V^{\u\d}$ 
has finite length.  Thus so does $V^\u$. In fact, the length of $V^{\u\d}$ 
is bounded by some function depending only on $d$ and $t$ in this situation.

It now follows from \ref{t-LO-FIP} that $H$ has $s$-LO for Prim, for some $s$ 
depending only on $t$ and $d$.
Using \ref{SpecPrim} we can show the same for Spec.  Let $Q\in\Spec(R)$ and
 let $I=(Q:H)$. Then $\widehat{Q}\in\Prim(\widehat{R})$. Thus  there exist 
$P_1, \dots ,P_n\in \Prim(\widehat{A})$ such that for each
$i$, $P_i\cap \widehat{R}=(\widehat{Q}:H)$ and $(\cap_i P_i)^t \subseteq
(\widehat{Q}:H)\widehat{A}$. Let $P'_i=P_i\cap A$. Then $P'_i\in\Spec(A)$
 and we have $P'_i\cap R=P_i\cap \widehat{R}\cap R=(\widehat{Q}:H)
\cap R=\widehat{I}\cap R=I$, and
$(\cap_i P'_i)^t \subseteq(\cap_i P_i)^t\cap A\subseteq (\widehat{Q}:H)
\widehat{A}\cap A\subseteq \widehat{I}\widehat{A}\cap A\subseteq
\widehat{IA}\cap A=IA$. Thus $H$ has $s$-LO for Spec. This answers an open problem 
from \cite{montgomery-schneider;prime}.

 Since $H$ also has coINC, GU, coGU and  coLO for Spec by
 the above, in order to show that all Krull relations hold it will be necessary
 only to prove INC.

 Now suppose that $R$ is an $H$-prime $H$-module algebra.  Then combining the 
above with \ref{H0prime} we have  

\begin{itemize}
\item $Q\in\Spec(R)$ is minimal if and only if $(Q:H)=0$, 
\item the minimal primes of $R$ are all conjugate under $G(H)$, there
 are at most $d$ of them, and their intersection is nilpotent of index at most 
$t+1$, and
\item there is a bijection between $H\Spec(R)$ and the space of orbits of $G(H)$
 on $\Spec(R)$.
\end{itemize}

More generally, we may obtain all the results above for {\em virtually pointed} $H$
 (that is, those which become pointed after a finite field extension). This
 includes all cocommutative $H$.

If both $H$ and $H^*$ are (virtually) pointed then $H$ satisfies all Krull relations 
for both Prim and Spec, as well as the finite induction and restriction properties.

\subsection*{Connected Hopf algebras} Suppose that $H$ is connected. Then since coINC is
satisfied, it follows from \ref{H0prime} that the map $Q\mapsto (Q:H)$ is a poset isomorphism
between $\Spec(R)$ and $H\Spec(R)$. Furthermore this map respects primitive and maximal ideals.
By the pointed case above, if $R$ is semiprime
 and $H$-prime then $R$ is in fact prime. This generalizes \cite[Lemma
1.2(ii)]{bergen-montgomery-passman}. 
 Dually, if $H^*$ is connected then $P\mapsto P\cap R$ is a poset isomorphism 
between $\Spec(A)$ and $H\Spec(R)$.  If also the extension is centralizing then we 
obtain an isomorphism between $\Spec(A)$ and $\Spec(R)$. In  particular, if
 $A$ is semiprime then so is $R$. Note that though this map is given
 by contraction,  its inverse is not given by expansion. Indeed, $0$
 is not a prime ideal for any nontrivial Hopf algebra, yet the
 extension $F\subseteq H$ is centralizing. Applying the correspondence
 above to this latter extension recovers the easy fact that $H$ must
 be a local ring with unique maximal ideal its augmentation ideal.

If $H$ and $H^*$ are both connected then again $\Spec(R)$ and $\Spec(A)$ are isomorphic, and all 
Krull relations and the finite induction and restriction properties are satisfied. We present a
simple example of such an $H$. Let
$H$ be the truncated divided power Hopf algebra in one variable over $F$.  As an algebra, $H$ is
spanned by elements $x^{(i)}$ for $0\leq i\leq p-1$. Here $x^{(0)}=1$ and the multiplication is
given by $x^{(i)}x^{(j)}=\binom{i+j}{j}x^{(i+j)}$. Thus as an algebra, $H\cong F\langle t\mid
t^p=0\rangle$. The comultiplication is given by making $t=x^{(1)}$ primitive, the augmentation
ideal is precisely $(t)$ and the antipode is given by $S(t)=-t$. 

Now $H$ is connected and self-dual. A comodule algebra
$A$ for $H$ is then just an algebra equipped with a derivation $d$ which is nilpotent of index
$p$. The coinvariants are the constants for $d$. The extension is faithfully flat Galois if and
only if there is an element $a\in A$ with $d^{p-1}(a)=1$. 

\subsection*{Semisimple Hopf algebras}

Montgomery and Witherspoon \cite{montgomery-witherspoon;irreducible-crossed}
 introduced the concept of {\em semisolvable} Hopf algebra, namely one  with  a subnormal series
all of whose factors are commutative or cocommutative. Classical structure theorems show that a
semisimple commutative Hopf algebra, after an
 algebraic extension $E$ of the ground field $F$, has the form $(EG)^*$. In characteristic zero,
after such an extension, a  cocommutative semisimple Hopf algebra has the form $EG$ (in
characteristic
 zero), while in characteristic $p>0$ it must be of the form $(EL)^*\#EG$,
 where $L$ is a $p$-group. In each case, the factors satisfy all Krull
 relations, as well as the  finite induction and restriction properties. In
 characteristic zero, the semisimple induction and restriction properties
 also hold for the factors. It now follows that every semisimple
 semisolvable Hopf algebra in characteristic zero is strongly semiprimitive
 and strongly cosemiprimitive, and has all Krull relations.

All known semisimple Hopf algebras of dimension less than 60, as well as many other examples, are
semisolvable \cite{montgomery;classifying-semisimple}. Thus the most pressing test question is
whether the Hopf algebras  constructed by Nikshych in
\cite{nikshych;twisting},  as deformations of the group algebra of the alternating group $A_5$,
 satisfy INC and the semisimple induction property.

\section{Conclusion}
\tlabel{sec:conc} In this paper it has been shown that the Krull relations for Prim imply those
for Spec, and that  strongly semiprimitive Hopf algebras are strongly semiprime. Thus in some sense
 it is superfluous to consider Spec and the whole theory can be founded on Prim. Since
strong semiprimitivity is implied by an appealing module property, it seems clear that the
behaviour of simple modules in faithfully flat Galois extensions should be a priority for
further work.

Many conjectures are suggested by the previous sections. Although many of  them are obvious and
it is highly unclear how to attack them, nevertheless
 we shall list some  explicitly (as questions) in the hope of stimulating
 further work.

 Let $H$ be a finite-dimensional Hopf algebra.
\begin{enumerate}
\romenumi
\item Does $H$ satisfy all Krull relations for Prim?
\item Does $H$ satisfy all Krull relations for Spec?
\item Does $H$ satisfy the finite induction (restriction) property? 
\end{enumerate} The answer to none of these questions is known, even for semisimple, pointed, or
connected Hopf algebras. Apart from group algebras and their duals, there is  no ``naturally
occurring" class of Hopf algebras for which all Krull relations are known to be satisfied. As
noted above, the missing relation for virtually pointed $H$ is INC.  In the
cocommutative case each Krull relation can be reduced (as in Section~\ref{sec:examples}) to the
special case of restricted enveloping algebras $u(L)$. However, all Krull relations hold when
$L$ is solvable (by transitivity and the results for abelian $L$ obtained by Chin
\cite{chin;prime-differential}). Thus a good test problem is to try to show that INC  holds for
$u(sl_2)$.  

We conclude with two problems on semisimple $H$. 
\begin{enumerate}
\item[(iv)] Does every semisimple $H$ satisfy the semisimple induction property? 
\item[(v)] Is every semisimple $H$ strongly semiprimitive?
\end{enumerate}

I thank Susan Montgomery for her hospitality at the University of Southern  California where
this work was started, and her and Hans-J{\"u}rgen Schneider for helpful discussions,
particularly on the proof of \ref{intcoGU}. I also thank the Department of Mathematics,
Statistics and Computer Science at the University of Illinois at Chicago for their hospitality
during the writing of this paper.
 
\bibliographystyle{amsalpha} 
\bibliography{hopf,primitive}  
  
\end{document}